\documentclass[12pt]{iopart}
\newtheorem{rmk}{Remark}[section]
\newtheorem{theo}{Theorem}[section]
\usepackage{amssymb}
\pdfoutput=1
\begin{document}

%\begin{frontmatter}

\title [Hypercomplexification of SDEs]{Exact integration of some systems of stochastic differential equations through hypercomplexification}
\author{C. Wafo Soh$^{1,2}$ and F. M. Mahomed$^2$}
\address{$^1$ Department of Mathematics  and Statistical Sciences, Jackson State University\\
 JSU Box 17610, 1400 JR Lynch Street, Jackson, MS 39217, USA
 \ead{Celestin.Wafo\_Soh@jsums.edu}} 
\address{$^2$DST-NRF Centre of Excellence in Mathematical and Statistical Sciences,\\ School of Computer Science and Applied Mathematics, University of the Witwatersrand, Johannesburg, Wits 2050, South Africa
\ead{Fazal.Mahomed@wits.ac.za}}
\begin{abstract}
We leverage commutative hypercomplex analysis to find closed-form solutions of  some systems of stochastic differential equations. Specifically, we obtain necessary and sufficient conditions under which a system of stochastic differential equations can be transformed into a scalar one involving processes valued in a commutative hypercomplex. In the event the targeted scalar stochastic differential equation is solved by quadratures, we recover the solution of the original system by projecting the solution of the scalar stochastic differential equation along the units of the underlying commutative hypercomplex. The conversion of a system of stochastic differential equations involving real-valued processes into a scalar one written in terms of hypercomplex-valued processes is termed {\em hypercomplexification}. Both hypercomplexification and its reverse are  mediated by the analyticity of  stochatic differential equations data. They may be iterated in order to generate higher-dimensional integrable systems of stochastic differential equations and  solve them. We showcase the utility of  hypercomplexification by treating several examples including linear, and linearizable systems of stochastic differential equations and stochastic Lotka-Volterra systems. Although we consider only random systems driven by white noises, hypercomplexification is fundamentally algebraic, and it readily extends to stochastic systems involving other types of disturbances.    
\end{abstract}

%\msc{60H10, 30G35, 13Pxx}

%\keywords{
%stochastic differential equations,
%numbers, hypercomplex analysis, generalized complex numbers, hypercomplexification, base equation}
%\submitto{\jpa}
\maketitle

\section{Introduction}
Many phenomena occur in environments governed by chance. The evolution of such phenomena is often described by stochastic differential equations (SDEs). We may roughly distinguish SDEs according to the type of noise that drives them. In this paper, we are primarily concerned with SDEs driven by the so-called white noise, which  are pervasive in physics, finance,
engineering and biology \cite{farm2000,blan2000,louz2001,huan2002,arat2003,spro2004,du2006,wije2008,zhu2009,come2012,pal2014}.  

Finding closed-form solutions of SDEs i.e. solutions expressible in terms of Lebesgue integrals and  stochastic integrals  of data is a formidable task. Even for scalar SDEs, there is not a general procedure for finding exact solutions comparable to their deterministic counterpart. Perhaps the largest class of scalar SDEs that is integrable by quadratures is comprised of scalar linear SDEs. Thus the strategy often adopted for finding exact solutions of a nonlinear scalar SDE is to ascertain whether it can be invertibly mapped to a linear SDE. Such an approach was pioneered by Gard \cite{gard} and extended by various researchers \cite{unal2007,unal2008,unal2010,unal2012,unal2013,mele2011}. 

Throughout this work, we shall consider systems of SDEs of the form 
\begin{equation}
dX_t  = a(t,X_t) dt + b(t,X_t) dW_t, \label{eq1:intro}
\end{equation}   
where $X_t(\omega) = X(t)(\omega) = X(t, \omega) = [X_{1:n}(t,\omega)]^T$, $a(t,X_t) = [a_{1:n}(t,X_t)]^T$,
$b(t,X_t) = [b_{ij}(t,X_t)]_{1\le i,j\le n}$, $W_t(\omega)= W(t,\omega)= [W_{1:n}(t,\omega)]^T$,  $\omega$ is the stochastic variable,  the $W_i$'s are independent  one-dimensional standard Wiener processes, and Eq. (\ref{eq1:intro}) is taken in It\^o's sense. Given an $n$-dimensional {\em commutative hypercomplex} (i.e. a finite-dimensional vector space equipped with a multiplication which is associative, commutative, distributive with respect to vector addition and compatible with the vector space field multiplication), with units $e_1, e_2, \ldots, e_n$, we want to find out the necessary and sufficient conditions on $a$ and $b$ such that under the substitution 
\begin{equation}
Z_t(\omega) = X_1(t,\omega)\, e_1 + X_2(t,\omega)\,e_2+\cdots + X_n(t,\omega)\, e_n, \label{eq2:intro}
\end{equation} 
the system (\ref{eq1:intro}) is transformed into an integrable scalar SDE of the form
\begin{equation}
d Z_t  = f(t, Z_t) dt + g(t, Z_t) d\mathbb{W}_t, \label{eq3:intro}
\end{equation}
with 
\begin{equation}
\mathbb{W}_t = W_1(t,\omega)\,e_1 + W_2(t,\omega)\,e_2 +\cdots + W_n(t,\omega)\, e_n. \label{eq4: intro}
\end{equation}
In the scheme outlined above, Eq. (\ref{eq3:intro}) will be called  the {\em base SDE} and the transformation of Eq. (\ref{eq1:intro}) into Eq. (\ref{eq3:intro}) via the replacement (\ref{eq2:intro}) will be termed {\em hypercomplexification}. Note that after solving the base SDE i.e. Eq. (\ref{eq3:intro}), the components of the solution of Eq. (\ref{eq1:intro}) are obtained by projecting  $Z_t$ along the units $e_1$ to $e_n$. 

The remainder of this paper is organized as follows. In Section 2, we introduce hypercomplex numbers and provide a succinct introduction to hyprcomplex analysis. Also, we show how hypecomplex analysis may be employed for seeking exact solutions of systems of SDEs. Section 3 is devoted to the uncovering of some systems of linear SDEs which can be solved in closed-form. Section 4 is dedicated to the problem of linearizing a system of SDEs through hypercomplexification. In Section 5, we explain how one can generate integrable stochastic Lotka-Volterra systems starting from  the scalar one. Finally, in Section 6, we summarize our findings and we elaborate on some problems we could not solve.  

\section{Hypercomplex analysis and hypercomplexification of stochastic differential equations}
Complex numbers and quaternions are particular instances of hypercomplexes. The general theory of hypercomplexes was developed and completed in the nineteen century. In this section, we provide a concise yet operational  account of the theory of hypercomplexes with an emphasis on commutative ones that are employed through this work. For an in-depth treatment of hypercomplexes, we refer the reader to the works 
\cite{pier1881,wedd1907,shaw1903,dick1914}.  

\subsection{Hypercomplex numbers}

A {\em hypercomplex}  $H$ is a finite-dimensional vector space equipped with a multiplication which is associative, distributive with respect to vector addition, compatible with the multiplication of the vector space field, and admits an identity that we will sometimes abusively denote $1$ when there is no risk of confusion. In this section, we shall denote the field of the vector space $\mathbb{F}$. In practice, we shall commonly take $\mathbb{F}$ as $\mathbb{R}$ or $\mathbb{C}$. When the multiplication of a hypercomplex is commutative, the hypercomplex is called a 
{\em commutative hypercomplex}. For all the applications envisaged in this paper, we solely use commutative hypercomplexes.

Now, suppose that  we have fixed a basis of a hypercomplex $H$, $\{e_1,e_2,\ldots, e_n\}$, so that 
\begin{equation}
H =\{ z = x_1\,e_1+x_2\,e_2 + \ldots + x_n\,e_n \; | \; (x_1, x_2, \ldots, x_n)\in \mathbb{F}^n \}. \label{eq1:hyp} 
\end{equation}
The elements of a basis of $H$ are called the {\em units} of the hypercomplex. Note that we can always select one of the units to be the identity of the hypercomplex. Owing to the closure property of the multiplication of $H$ which states that the product of two generic element of $H$ must belong to $H$, we must have $e_i\, e_j \in H$ for all $i,j = 1:n$. Thus,  
\begin{equation}
e_i \, e_j = \sum_{k=1}^n  \gamma_{ijk} e_k,\; \mbox{for all }i, j = 1:n,  \label{eq2:hyp}
\end{equation}
where the $\gamma_{ijk}$'s belong to $\mathbb{F}$. The $\gamma_{ijk}$'s are the {\em structure constants} of the hypercomplex $H$ associated with the units $e_1, e_2,\ldots, e_n$. It can be shown that the multiplication of $H$ is commutative if and only if 
\begin{equation}
\gamma_{ijk} = \gamma_{jik},\; \mbox{for all } i,j,k = 1:n. \label{eq3:hyp}
\end{equation}          
It can be verified that the multiplication of $H$ is associative if and only if $(e_ie_j) e_k = e_i (e_je_k),\; 
\mbox{ for all } i,j,k = 1:n$, which translates into 
\begin{equation}
\sum_{s=1}^n \gamma_{ijs} \gamma_{sk t} = \sum_{s=1}^n \gamma_{ist}\gamma_{jks}, \;\mbox{for all } i,j,k,t = 1:n. \label{eq4:hyp}
\end{equation}
The element $\epsilon = \sum_{i=1}^n \epsilon_i e_i$ is the multiplicative identity of $H$ if and only if $\epsilon e_i =
e_i$ and $e_i \epsilon = e_i$ for all $i = 1:n$. That is 
\begin{equation}
\sum_{k=1}^n  \epsilon_k \gamma_{ikj} = \sum_{k=1}^n \epsilon_k \gamma_{kij} = \delta_{ij} 
\mbox{ for all } i, j = 1:n,\label{eq4b:hyp}
\end{equation}
where $\delta_{ij}$ is Kronecker's symbol.
In practice, a hypercomplex is specified by providing its multiplication table i.e. a table in which the entry located on the $i^{\,\mbox{th}}$ row and  $j^{\,\mbox{th}}$ column, is $e_i e_j$. For commutative hypercomplexes, such a table is symmetric which respect to its main diagonal.

The classification of lower-dimensional hypercomplexes was done by several researchers 
\cite{pier1881,dick1914,sche1889,stud1890,hawk1902,pasq1920} using often different or complementary techniques. In our recent paper \cite{wafo2016}, we have summarized  from the work of Study \cite{stud1890} the classification of commutative hypercomplexes up to dimension four. We shall use the nomenclature of \cite{wafo2016} for lower-dimensional commutative hypercomplexes.

There are three nonequivalent classes of two-dimensional commutative hypecomplexes over $\mathbb{R}$ that can be treated at once through the consideration of the so-called generalized complex numbers \cite{yagl2014,hark2004}.  The set of general hypercomplex numbers is denoted $\mathbb{C}_p$ and its units may be chosen as $e_1 = 1$, $e_2 = i$ with its multiplication table given by:
$e_1^2 = 1,\; e_2^2 = p,\; e_1e_2 = e_2e_1 =e_2$, where $p$ is a real number. Put another way $\mathbb{C}_p =\{x+ iy \,| \,x, y\in \mathbb{R}\}$ and the algebra is done in $\mathbb{C}_p$ like in $\mathbb{R}$ under the proviso that $i^2 = p$. 

According to Study's \cite{stud1890} classification, there are up to equivalence transformations,  five  three-dimensional and twelve four-dimesional commutative hypercomplexes over $\mathbb{R}$. Note that some high-dimensional  hypercomplexes  may be generated through direct sum or direct product \cite{dick1914, wafo2016} of lower-dimensional ones.
\subsection{ A primer on hypercomplex analysis} 
The investigation of analytic functions of hypercomplex variables was initiated by Scheffers \cite{sche1893,sche1893b} for commutative hypercomplexes. His works was extended in several directions by Ketchum \cite{ketc1928} who in particular generalized all important theorems of the classical complex analysis to  functions of hypercomplex variables in commutative hypercomplexes. Our aim in this subsection is to provide a brief account of this theory. We shall restrict ourselves to situations that suit our need in this paper. For an in-dept treatment of commutative hypercomplex analysis, we refer the reader to the seminal papers \cite{sche1893,sche1893b,ketc1928}.

Throughout this subsection, $H$ will be a commutative hypercomplex over real numbers, with units $e_1, e_2,\ldots, e_n$, and corresponding structure constants $\gamma_{ijk},\; i,j,k =1:n$. We denote the multiplicative identity of $H$ by $\epsilon = \sum_{i=1}^n \epsilon_i \, e_i$. We topologize $H$ trough the Euclidian norm defined for $z= \sum_{i=1}^{n} x_i \, e_i$ by 
$||z|| = \sqrt{\sum_{i=1}^n x_i^2}$.

A function of a hypercomplex variable  $f(z) = \sum_{i=1}^n f_i (x_{1:n})\, e_i$ is differentiable at $z_0 \in H$ if there is $u\in H$ such that 
\begin{equation}
\lim_{||h||\rightarrow 0}  ||f(z +h) - f(z) - u\,h|| =0. \label{eq1:han} 
\end{equation}
We call $u$ the derivative of $f$ at $z_0$ and we denote it by $f'(z_0)$.  

The function $f$ is {\em monogenic} in a region (open connected set)  $U$ of $H$ if it is differentiable at every point in U and the $f_i$ are analytic on $U$. 

Now assume that $f$ is monogenic in a region $U$. Then, on this region, we must have
\begin{equation}
df(z) = f'(z) dz. \label{eq2:han} 
\end{equation}       
Substitute  $f'(z) = \sum_{i=1}^n f'_i(z)e_i$, $dz = \sum_{i=1}^n dx_i e_i$ into (\ref{eq2:han}) to obtain 
\begin{equation}
\sum_{j=1}^n {\partial f_i \over \partial x_j} dx_j = \sum_{j,k} \gamma_{kji} f'_k \,dx_j. 
\label{eq3:han}
\end{equation}
Since Eq. (\ref{eq3:han}) must hold true for all increment $dx_j$, $j= 1:n$, we must have
\begin{equation}
{\partial f_i \over \partial x_k} = \sum_{j=1}^n \gamma_{jki}f'_j \mbox{ for all } i,k = 1: n.
\label{eq4:han}
\end{equation}
Multiply both sides of Eq. (\ref{eq4:han}) by $\epsilon_k$ (the $k$th component of the multiplicative identity) and sum both sides from $k=1$ to $k=1$. Then employ the identity (\ref{eq4b:hyp}) on the left-hand side to arrive at  
\begin{equation}
f'_i = \sum_{k=1}^n \epsilon_k\, {\partial f_i \over \partial x_k} \mbox{ for all } i =1:n.  \label{eq5:han}
\end{equation}
Substituting  Eq. (\ref{eq5:han}) into Eq. (\ref{eq4:han}) yields
\begin{equation}
{\partial f_i \over \partial x_k} = \sum_{j,\ell = 1}^n \epsilon_{\ell} \gamma_{jki}{\partial f_k \over x_{\ell}}\,\cdot
\label{eq6:han}
\end{equation} 
Note that the system (\ref{eq6:han}) generalizes the classical Cauchy-Riemann equations. We shall call it {\em Scheffers' equations} owing to the following theorem proved by Scheffers \cite{sche1893,sche1893b}.
\begin{theo}
The function $f(z) = \sum_{i=1}^n f_i(x_{1:n})\, e_i$, where $z= \sum_{i=1}^n x_i e_i \in H$, and $H$ is a commutative hypercomplex with structure constants $\gamma_{ijk},\; i,j,k = 1:n$, is monogenic on a region $U$ if and only if the $f_i$'s are continuously infinitely differentiable on $U$  and the system (\ref{eq6:han}) is satisfied on $U$. 
\end{theo}

A function of a hypercomplex variable $z$ is {\em analytic} at a point $z_0$ if it can be expressed as the sum of a power series in $(z-z_0)$ in the vicinity of $z_0$. Such a function is analytic in a region $U$ if it is analytic at every point in $U$. It can be verified that analytic functions are monogenic. In fact, Scheffers \cite{sche1893,sche1893b} proved that monogenic functions are analytic. Thus, in the sequel, we shall use monogenic and analytic as synonymous when dealing with functions of  hypercomplex variables.

Several classical analytic functions have their hypercomplex counterparts which are also analytic. The machinery commonly employed to achieve the hypercomplexification of these functions  includes series expansion, analytic continuation or integral representation.  In what follows, we shall illustrate this point by constructing the exponential functions on $\mathbb{C}_p$  and $A_3^4$ \cite{wafo2016} through series expansion.

Let $z = x + i y \in \mathbb{C}_p$. We define the exponential of $z$ as follows \cite{wafo2016,yagl2014,hark2004}:
\begin{eqnarray}
\mbox{e}^z &=& \mbox{e}^x \mbox{e}^{i y} \nonumber\\
&=& \mbox{e}^x \sum_{k =1}^{\infty} {(iy)^k\over k!} \nonumber\\
& = & \mbox{e}^x \left \{  \sum_{k =1}^{\infty} {(iy)^{2k}\over (2k)!}  +
 \sum_{k =1}^{\infty} {(iy)^{2k+1}\over (2k+1)!}  \right \} \nonumber \\
 & = &  \mbox{e}^x \left \{  \sum_{k =1}^{\infty} {p^k y^{2k}\over (2k)!}  +
 i\sum_{k =1}^{\infty} {p^k y^{2k+1}\over (2k+1)!}  \right \} \nonumber \quad [ i^2 = p] \\
& = & \mbox{e}^x [\cos_p (y) + i \sin_p (y) ], \label{eq7:han}  
\end{eqnarray}    
where we have set
\begin{equation}
\cos_p(y) = \sum_{k =1}^{\infty} {p^k y^{2k}\over (2k)!}     
\mbox{ and } \sin_p(y) = \sum_{k =1}^{\infty} {p^k y^{2k+1}\over (2k+1)!} \; \cdot \label{eq8:han}
\end{equation}
By noticing that $\mbox{e}^{iy}\mbox{e}^{-iy} =1$, we arrive at the generalized Pythagorean identity
\begin{equation}
\cos_p(y) - p\, \sin_p(y) = 1. \label{eq9:han} 
\end{equation}

Let us now turn our attention to the commutative hypercomplex  $A^4_3$ \cite{wafo2016} whose units are $e_1 =1$, $e_2=i$ and $e_3 = j$ and the multiplication table is given by: $i^2=j$, $ij = 0$ and $j^2 =0$. Let us express a generic element of $A_3^4$ as $z = t + x i + y j$. So, we may define $\mbox{e}^z$ as
\begin{eqnarray}
\mbox{e}^z &=& \mbox{e}^t \mbox{e}^{ix} \mbox{e}^{jy}\nonumber \\
& = & \mbox{e}^t \; \sum_{k=1}^{\infty} {(ix)^k\over k!}  \sum_{k=1}^{\infty}\; {(jx)^k\over k!}\nonumber \\
&= & \mbox{e}^t \; \left (1+ x\,i + {x^2\over 2} j \right )(1+ j y)\nonumber \\
&=& \mbox{e}^t\left [ 1 + x \, i + \left ({x^2\over 2} +y \right ) j\right ]. \label{eq10:han}  
\end{eqnarray}
Then, the natural logarithm, $\ln z$, is obtained  by solving the equation  $z = \exp (a + b\, i + c\, j)$ for $a$, $b$ and $c$. The foregoing equation yields the system 
\begin{eqnarray}
 t & = & \mbox{e}^a , \label{eq11:han}\\
 x & = & \mbox{e}^a b, \label{eq12:han}\\
 y & = & \mbox{e}^a \left ({b^2 \over 2} + c \right ). \label{eq13:han}  
\end{eqnarray}
The solution of the system (\ref{eq11:han})-(\ref{eq13:han}) is given by: $a  =\ln t$, $b = x/t$, and 
$c = y/t - x^2/ (2t^2)$. Thus, the natural logarithm on $A_3^4$ is defined as follows:
\begin{equation}
\ln z = \ln t + {x\over t}\, i + \left ( {y\over t} -{x^2\over 2 t^2}\right ) j. \label{eq14:han}
\end{equation}
From Eqs. (\ref{eq10:han}) and (\ref{eq14:han}), we infer that for any real number $m$,
\begin{equation}
z^m = \exp(m \ln z)= t^m \left [ 1 + {m \,x \over t} i + \left ( {m\, y \over t} + {m(m-1) x^2 \over 2 t^2} \right ) j \right ] . \label{eq15:han}
\end{equation}        
In fact we can even define $z_1^{z_2}$ for $z_1, z_2 \in A^4_3$ by noticing that $z_1^{z_2} = \exp(z_2\ln z_1)$ and using the formulas (\ref{eq10:han}) and (\ref{eq14:han}).

From the formula (\ref{eq15:han}), we deduce that 
\begin{eqnarray}
\cos(z) & = & \sum_{m=0}^{\infty} (-1)^m {z^{2m}\over  (2m)!} = \sum_{m=0}^{\infty} (-1)^m\, {t^{2m}\over (2m)! }
+ i\; \sum_{m =1}^{\infty} (-1)^m\, {t^{2m-1} x \over (2m-1)!}\nonumber \\
& & + j \;\sum_{m=0}^{\infty} {(-1)^m \,\over (2m)!}\left (2m \,y\, t^{2m-1} + 2 m(2m-1)\, t^{2m-2}\,{x^2\over 2}\right ) \nonumber\\
& = & \cos t -  (x \sin t)i  -\left ( y \sin t +{x^2 \cos t\over 2}\right ) \, j. \label{eq16:han} 
\end{eqnarray}
Likewise, we find that 
\begin{equation}
\sin(z) =  \sin t +  (x \cos t)\,i  +\left ( y \cos t -{x^2\sin t\over 2}\right ) \, j. \label{eq17:han}
\end{equation}
The reader may verify that amazingly $\cos^2 z + \sin^2 z = 1$ for all $z\in A_3^4$ such that the latter equation makes sense. 
 
\subsection{Hypercomplexification of stochastic differential equations} 
It is hard to find the general exact solution of a scalar SDE. Usually, exact solutions of scalar ODEs are found through clever use of the stochastic chain rule (It\^o's lemma) to map these SDEs to integrable SDEs. Arguably, the largest class of scalar integrable SDEs is the class made of linear SDEs. Thus the main strategy for finding exact solutions of a scalar SDE consists in looking for invertible change of both the dependent and dependent variables that convert the underlying SDE into a linear one. Such an approach was pioneer by Gard \cite{gard} and extended by several researchers
\cite{unal2007,unal2008,unal2010,unal2012,unal2013,mele2011}. 

For systems of SDEs, the situation is even more complicated than for scalar SDE. Indeed one has to deal with the coupling of the equations and the fact that in general, linear system of SDEs are not exactly integrable. So even if one contemplates transforming a system of SDE to a linear system of SDEs, the target system may not be exactly integrable.

Our proposal for the search of integrable systems of SDEs consists in ascertaining whether by  appropriate hypercomplex substitutions one can transforms these systems to scalar SDEs in which the involved process are valued in  commutative hypercomplexes. In the event one can transform a system of SDEs to a scalar SDE and the scalar SDE is solvable by quadratures, one recovers the solution of the original system by projecting the solution of the target scalar SDE along the units of the commutative hypercomplex: The exact solution of the scalar SDE will be expressed in terms of Lebesgue and stochastic integrals involving data and elementary functions whose hypercomplex counterpart can be generated as discussed in the previous subsection. 
In this scheme, we shall call the target scalar SDE the {\em base equation} \cite{qadi2015}. The conversion of a systems of SDEs in which the dependent variables and data are real-valued to a scalar SDE which has hypercomplex-valued dependent variable and data, will be termed {\em hypercomplexification}.

In order to reify the above narrative, consider the following scalar SDE
\begin{equation}
dZ(t) = a(t, Z(t))dt + b(t,Z(t))d\mathbb{W}(t), \label{eq1:xn} 
\end{equation}  
where $Z(t) =\sum_{i=1}^n X_i(t)\,e_i$ is valued in a commutative  hypercomplex over reals, $H$, with units $e_1, e_2,\ldots, e_n$,  
$\mathbb{W}(t) = \sum_{i=1}^m W_i(t)\,e_i$, and $m\le n$. In other sections, we shall take $m=n$ for the sake of simplicity. By letting $a = \sum_{i=1}^m a_i(t,X_{1:n}(t))\,e_i$ and   
$b = \sum_{i=1}^m b_i(t,X_{1:n}(t))\,e_i$,  and using the linear independent of the units of $H$, the equation (\ref{eq1:xn}) implies
\begin{equation}
dX_i(t)  =  a_i(t, X_{1:n}(t))dt  
 +\sum_{k=1}^n \sum_{j=1}^m \gamma_{ijk}b_j(t, X_{1:n}(t))dW_k(t).
\label{eq2:xn}
\end{equation}
The scalar equation (\ref{eq1:xn}) is the {\em base equation}  of the system (\ref{eq2:xn}). 
In order to reverse the procedure i.e. start from the system (\ref{eq2:xn}) and arrive at Eq. (\ref{eq1:xn}), one must require that the $\gamma_{ijk}$'s be the structure constants of a commutative hypercomplex and  that the functions $a$ and $b$ be analytic in $Z$ as functions valued in that hypercomplex. Thus, these functions must satisfy Scheffers' equations (\ref{eq6:han}).

In practice, one generates Eq. (\ref{eq2:xn}) by using the classification of lower-dimensional commutative hypercomplexes \cite{stud1890,wafo2016} and a base equation one knows how to solve. We shall illustrate this procedure in the remaining sections.

Note that the above procedure may be iterated in such a way that one may end up with a high-dimensional integrable
system of SDEs. Indeed, if the system (\ref{eq2:xn}) is integrable, we can also hypercomplexify it using another hypercomplex, $K$, say. The resulting system could be generated from the base equation (\ref{eq1:xn}) using the direct product \cite{wafo2016} of $H$ by $K$, $H\otimes K$.  Thus stating from a scalar integrable SDE, we may build a hierarchy of integrable systems of SDEs that may be nontrivial to solve without employing the hypercomplexification procedure.

\section{Integrable systems of linear stochastic differential equations}
In this section, we characterize linear systems of SDEs that are integrable by hypercomplexification. It is important to stress  that  although scalar linear SDEs are integrable, systems of linear SDEs do not in general admit closed-form  solutions \cite{gard}.
Consider a scalar linear SDE
\begin{equation}
dZ(t) = (f_1(t)+ f_2(t)Z(t))dt + [g_1(t)+ g_2(t)Z(t)] d\mathbb{W}(t) \label{eq1:lsde} 
\end{equation}
in which the functions $f_1$, $f_2$, $g_1$ and $g_2$ are deterministic functions and we have dropped the stochastic variable as it is customary in stochastic calculus. 
Its closed-form solution is \cite{gard}
\begin{eqnarray}
Z(t) &=& E(t)\left \{  Z(0) + \int_0^t E^{-1}(s) [ f_1(s) - g_1(s)g_2(s)]ds \right. \nonumber \\
& & \left. + \int_0^t E^{-1}(s)g_1(s)d\mathbb{W}(s) \right \},
\label{eq2:lsde}
\end{eqnarray}
where
\begin{equation}
E(t) = \exp \left \{ \int_0^t \left [ f_2(s) -{1\over 2} g_2^2(s)\right ]ds + \int_0^t g_2(s)d\mathbb{W}(s) \right \}.
\label{eq3:lsde}
\end{equation}
Now, let us treat the scalar quantities as elements of a commutative hypercomplex $H$ with units $e_1$ to $e_n$, and structure constants $\gamma_{ij}$ where $i,j = 1:n$:
\begin{eqnarray}
Z(t) & = & \sum_{i=1}^n X_i (t)e_i, \label{eq4:lsde}\\
f_1 (t)  & = & \sum_{i=1}^n f_{1i}(t)\,e_i, \label{eq5:lsde}\\
f_2 (t)  & = & \sum_{i=1}^n f_{2i}(t)\,e_i, \label{eq6:lsde}\\
g_1 (t)  & = & \sum_{i=1}^n g_{1i}(t)\,e_i, \label{eq7:lsde}\\
g_2 (t)  & = & \sum_{i=1}^n g_{2i}(t)\,e_i, \label{eq8:lsde}\\
d\mathbb{W}(t) & = & \sum_{i=1}^n dW_i(t) e_i. \label{eq9:lsde}
\end{eqnarray}
Now, substitute Eqs. (\ref{eq4:lsde})-(\ref{eq9:lsde}) into Eq. (\ref{eq1:lsde}) to obtain after some algebraic manipulations
\begin{eqnarray}
dX_i(t) & =&  \left (f_{1i}(t) + \sum_{k,\ell =1}^n \gamma_{k\ell i}\, f_{2k}(t) X_{\ell}(t)  \right )dt \nonumber\\
& & + 
\sum_{k,\ell=1}^n \left (\gamma_{k\ell i}\, g_{1k}(t)   + \sum_{m, p =1}^n \gamma_{kp m} \gamma_{m\ell i} \,g_{2k}(t)X_{p}(t)  \right ) dW_{\ell}(t), \label{eq10:lsde}
\end{eqnarray} 
for all $i = 1:n$. Conversely, starting from the system (\ref{eq10:lsde}) and provided that the $\gamma_{ijk}$'s appearing in it are structure constants of a commutative hypercomplex with basis $e_1, e_2,\dots, e_n$, we can transform Eq. (\ref{eq10:lsde}) into Eq.(\ref{eq1:lsde}) via the substitution (\ref{eq4:lsde}). We has thus established the following theorem.
\begin{theo}
\label{theo1:lsde}
A system of $n$ linear SDEs that can be hypercomplexified provided  has to have the form (\ref{eq10:lsde}), where the $\gamma_{ij}^k$'s are structure constants of a commutative hypercomplex.  The solution of the system (\ref{eq10:lsde}) is obtained by projecting Eq.(\ref{eq2:lsde}) along the units of the underlying commutative hypercomplex.
\end{theo}
\begin{rmk} Stochastic differential equations of the form
\begin{equation}
du_i = - \sum_{j=1}^3 A_{ij}(t) u_{j} dt + \sum_{j=1}^3 B_{ij}(t) dW_j,\; i= 1:3 \label{eq10b:lsde}
\end{equation} 
arise in the study of the dispersion of particles in turbulent flows \cite{pope2002}. By employing three-dimensional hypercomplexes, we may infer from the system (\ref{eq10:lsde}) integrable forms of Eq. (\ref{eq10b:lsde}). 
\end{rmk}

Let us illustrate Theorem \ref{theo1:lsde} by using $\mathbb{C}_p$ as our commutative hypercomplex. For $\mathbb{C}_p$, the system (\ref{eq10:lsde}) becomes:
\begin{eqnarray}
dX_1(t) & = & \left [ f_{11} (t)+ f_{21}(t) X_1(t) + p\, f_{22}(t)X_2(t) \right ] dt \nonumber \\
& & +\left [ g_{11}(t) + g_{21}(t)X_1(t) + p\, g_{22}(t) X_2(t) \right ] dW_1(t)\nonumber \\
& & + p\,\left [ g_{12}(t)+g_{22}(t) X_1(t) + g_{21}(t)X_2(t)\right ] dW_2(t), \label{eq11:lsde}\\
dX_2(t) & = & \left [ f_{12}(t) + f_{22}(t) X_1(t) + f_{21}(t) X_2(t) \right ] dt \nonumber \\
& & +\left [ g_{12}(t)  + g_{22}(t)X_1(t) + g_{21}(t)X_2(t)\right ] dW_1(t)\nonumber \\
& & + \left [g_{11}(t)+ g_{21}(t)X_1(t) + p\, g_{22}(t) X_2(t)\right ] dW_2(t). \label{eq12:lsde}
\end{eqnarray}
In order to express the solution of the system (\ref{eq11:lsde})-(\ref{eq12:lsde}) in closed-form using Eq.(\ref{eq2:lsde}), we introduce the following notations
\begin{eqnarray}
a(t,\tau) & = & \int_{\tau}^t \left [ f_{21}(s) -{1\over 2} (g_{21}^2(s) + p g_{22}^2(s)) \right ] ds\nonumber \\
& & + \int_{\tau}^t g_{21}(s)dW_1(s) + p \int_{\tau}^t g_{22}(s) dW_2(s), \label{eq13:lsde}\\
b(t,\tau) & = & \int_{\tau}^t [ f_{22}(s) - g_{21}(s) g_{22}(s)]\,ds \nonumber \\
& & + \int_{\tau}^t g_{22}(s) dW_1(s) + \int_{\tau}^t g_{21}(s)dW_2(s).
\label{eq14:lsde}
\end{eqnarray}
Then, the general solution of the system (\ref{eq11:lsde})-(\ref{eq12:lsde}) is given by
\begin{eqnarray}
 X_1(t)  &=&  \mbox{e}^{a(t,0)}[X_1(0)\cos_p[b(t,0)] + p X_2(0)\sin_p[b(t,0)] \nonumber \\
& & + 
\int_0^t  \mbox{e}^{a(t,s)} \left [   \cos_p[b(t,s)] \{f_{11}(s) -g_{11}(s)g_{21}(s) -pg_{12}(s)g_{22}(s)\} 
\right.   \nonumber  \\
& & \left. + p \sin_p[b(t,s)]\{f_{21}(s) -g_{12}(s) g_{21}(s) -g_{11}(s)g_{22}(s)\}\right ]  ds 
\nonumber   \\
& & + \int_0^t \mbox{e}^{a(t,s)} \left [\cos_p [b(t,s)] g_{11}(s) + p \sin_p[b(t,s)]g_{12}(s)\right ] dW_1(s) \nonumber\\
& & + p\int_0^t \mbox{e}^{a(t,s)} \left [\cos_p [b(t,s)] g_{12}(s) +  \sin_p[b(t,s)]g_{11}(s)\right ] dW_2(s)
\label{eq15:lsde},\\
X_2(t) & = & \mbox{e}^{a(t,0)}[ X_1(0)\sin_p[b(t,0)] + X_2(0)\cos_p[b(t,0)]] \nonumber \\
& & + 
\int_0^t \left \{\mbox{e}^{a(t,s)} \left [   \sin_p[b(t,s)] \,\{f_{11}(s) -g_{11}(s)g_{21}(s) -pg_{12}(s)g_{22}(s)\} \right. \right.  \nonumber  \\
& & \left. \left. +  \cos_p[b(t,s)]\,\{f_{21}(s) -g_{12}(s) g_{21}(s) -g_{11}(s)g_{22}(s)\}\right ] \right \} ds \nonumber\\
& & + \int_0^t \mbox{e}^{a(t,s)} \left [\cos_p [b(t,s)]\, g_{12}(s) +  \sin_p[b(t,s)]\,g_{11}(s)\right ] dW_1(s) \nonumber\\
& & + \int_0^t \mbox{e}^{a(t,s)} \left [\cos_p [b(t,s)] \,g_{11}(s) +  p\sin_p[b(t,s)]\,g_{12}(s)\right ] dW_2(s).
\label{eq16:lsde}
\end{eqnarray}
\begin{rmk} (1) It can be verified that the system (\ref{eq11:lsde})-(\ref{eq12:lsde}) generalizes both the circular and hyperbolic Brownian motions.

(2) We may infer from Eqs. (\ref{eq15:lsde})-(\ref{eq16:lsde}) that the fundamental matrix of the deterministic system 
$$ {dX_1\over dt} = f_{21}(t) X_1 + p\, f_{22}(t) X_2,\quad {d X_2\over dt} = f_{22}(t) X_1 + f_{21}(t) X_2$$
is 
$$ \Phi(t) = \exp\left (\int_0^t f_{21} (s) ds \right ) \left [
\begin{array}{ll}
\cos_{p} \left ( \int_0^t f_{22}(s) ds  \right ) & \quad p\sin_{p} \left ( \int_0^t f_{22}(s) ds  \right) \\
\sin_{p} \left ( \int_0^t f_{22}(s) ds  \right ) &\quad \cos_{p} \left ( \int_0^t f_{22}(s) ds  \right )
\end{array} \right ].
$$ 
It can be demonstrated  that $\det(\Phi(t)) = \exp\left (2\int_0^t f_{21} (s) ds\right )$ as one expects.
\end{rmk}

\section{Reducibility by hypercomplexification of  systems of   stochastic differential equations}
Here we consider the problem of hypercomplexifying a system of SDEs to a {\em reducible} scalar SDE. We recall that a scalar SDE $dZ(t) = f(t,Z(t)) dt  + g(t, Z(t))d\mathbb{W}(t)$ is reducible if it can be transformed to the linear scalar SDE
$dY(t) = a(t) dt + b(t)d\mathbb{W}(t)$ by means of a substitution $Y(t) = h(t, Z(t))$. Reducibility is thus a 
narrow-sense linearization. It was considered by Gard \cite{gard}. The general problem of linearization of a scalar SDE was completely solved by Meleshko \cite{mele2011}. Although we focus here on reducibility through hypercomplexification, all the results will apply to linearization mutatis mutandis. We start by recalling Gard's reducibility theorem.
\begin{theo}
The scalar SDE 
\begin{equation}
dZ(t) = f(t,Z(t)) dt  + g(t, Z(t))d\mathbb{W}(t)\label{eq1:red}
\end{equation}
can be invertibly mapped via the substitution $Y(t) = h(t, Z(t))$ to the scalar linear SDE
\begin{equation}
dY(t) = a(t) dt + b(t)d\mathbb{W}(t) \label{eq2:red}
\end{equation} if and only if $\partial N/ \partial Z = 0$, where
\begin{equation}
N = g\left \{ {1\over g^2}{\partial g \over \partial t} - {\partial \over \partial Z} \left (   {f\over g} \right )
+{1\over 2}\,{\partial^2 g \over \partial Z^2} \right \}. \label{eq3:red}
\end{equation} 
The functions $a$, $h$ and $b$ are then respectively given by
\begin{eqnarray}
a &=& \exp \left ( \int^t N(s)ds \right), \label{eq4:red}\\
{\partial h \over \partial Z} &=& {a\over g}, \label{eq5:red}\\
b  &=& {\partial h \over \partial t} + f {\partial h \over \partial Z} + {1\over 2} g^2 {\partial^2 h \over \partial Z^2}
\,\cdot \label{eq6:redd}
\end{eqnarray}   
\end{theo}
Now, we wish to find the necessary and sufficient conditions under which a system of SDEs can be transformed to a scalar reducible SDE through the substitution $Z = X_1e_1+X_2e_2 + \cdots + X_n e_n$, where  $Z$ belong to a commutative hypercomplex. Owing to the fact that the approach used does not depend on the dimension of the hypercomplex, we shall without loss of generality use the hypercomplex $\mathbb{C}_p$ so that the substitution will take the form 
$Z = X + i Y \in \mathbb{C}_p$. Starting from Eq. (\ref{eq1:red}), it can be shown through algebraic computations that a system of two SDEs with dependent variables $X$ and $Y$ can be hypercomplexified through the latter substitution if and only if  it assumes the form 
\begin{eqnarray}
dX(t) &=& f_1(t)\,dt + g_1(t, X(t), Y(t))\,dW_1(t) \nonumber \\
 & & + p g_2(t,X(t), Y(t))\, dW_2(t), \label{eq6:red}\\
dY(t) &=& f_2(t) \, dt  + g_2(t, X(t), Y(t))\,dW_1(t)\nonumber \\
& &  + g_1(t,X(t), Y(t))\,dW_2(t), \label{eq7:red}
\end{eqnarray}
with the $f_i$'s and $g_i$'s constrained as follows: 
\begin{eqnarray}
{\partial f_1 \over \partial X} &=& {\partial f_2\over \partial Y}, \quad 
{\partial f_2 \over \partial X} = p\, {\partial f_1 \over \partial Y}\\ \label{eq8:red}\\
{\partial g_1 \over \partial X} &=& {\partial g_2\over \partial Y}, \quad 
{\partial g_2 \over \partial X} = p\, {\partial g_1 \over \partial Y} \label{eq9:red}
\end{eqnarray}
The equations (\ref{eq8:red})-(\ref{eq9:red}) are Scheffers' equations for the hypercomplex $\mathbb{C}_p$. They ensure that the functions $f= f_1 + i f_2$ and $g = g_1 + i g_2$ are analytic so that the system (\ref{eq8:red})-(\ref{eq9:red})
can be converted back to a scalar SDE of the form (\ref{eq1:red}).
Now, since $f$ and $g$ are analytic thanks to (\ref{eq8:red})-(\ref{eq9:red}), all the functions appearing in Eq. (\ref{eq3:red}) are analytic. We may then split the foregoing equation as 
\begin{eqnarray}
N_1 &=& -f_{1,X} +{1\over 2}\, g_1 g_{1,XX} + {1\over 2}\, p\, g_2 g_{2,XX} \nonumber  \\
& & +\;{1\over g_1^2 -pg_2^2}  ( f_1g_{1,X} -p g_2 g_{1,X} + p f_2 g_{2,X} \nonumber \\
& &  - p f_1g_2g_{2,X}+g_1g_{1,t} -p g_2g_{2,X}  ) \label{eq10:red}\\
N_2 & = & -f_{2,X} +{1\over 2}\, g_2 g_{1,XX} +{1\over 2}\, g_1 g_{2,XX} \nonumber \\
& & + {1\over g_1^2 -p g_2^2} ( f_2g_1g_{1,X} -f_1g_2g_{1,X} + f_2 g_2 g_{2,X} \nonumber \\
& & + f_1 g_1 g_{2,X} -g_2 g_{1,X} + g_1 g_{2,t} ) \label{eq11:red}
\end{eqnarray}
Note that Eq. (\ref{eq4:red}) constraints $N$ to depend solely on $t$. Since the conditions on the $f_i$'s and $g_i$'s guarantee that $N$ is analytic,  the condition $\partial N/\partial Z = 0$  is equivalent to $\partial N_1 /\partial X = 0$ and $\partial N_2 / \partial X = 0$. However the latter conditions together with the Sheffers's equations lead to $N_1 = N_1(t)$ and $N_2 = N_2(t)$ only when $p\ne 0$. When $p= 0$, one has to impose additionally that $\partial N_{1}/ \partial Y = 0$ in order to guarantee that $N_1 = N_2(t)$.  Thus, we have proved the following theorem.
\begin{theo}
The system of SDEs (\ref{eq6:red})-(\ref{eq7:red}) can be transformed through the substitution $Z(t)= X(t) + i Y(t) \in \mathbb{C}_p$ to a scalar reducible SDE if and only if one of the following two conditions is satisfied.

(1) The real number $p$ is nonzero, the Eqs. (\ref{eq8:red})-(\ref{eq9:red}) are fulfilled, $\partial N_{1}/\partial X = 0$ and 
$\partial N_2 / \partial X = 0$.

(2) The number $p$ is zero,   the Eqs. (\ref{eq8:red})-(\ref{eq9:red}) are fulfilled, $\partial N_{1}/\partial X = 0$, 
$\partial N_2 / \partial X = 0$ and $\partial N_1/\partial Y = 0$. 
  
\end{theo}    

\section{Integrable stochastic Lotka-Volterra systems}
Stochastic  Lokta-Volterra systems are pervasive in the modeling of the competitive interaction of agents evolving in a random environment \cite{farm2000,blan2000,louz2001,huan2002,arat2003, spro2004,du2006,wije2008, zhu2009,come2012,pal2014}. The purpose of this section is to provide some stochastic Lokta-Volterra systems which can be integrated by quadrature through hypercomplexification. The notations and definitions that are not introduced here are those of the previous sections.

The scalar stochastic Lotka-Volterra model may be expressed as 
\begin{equation}
dZ(t) = (bZ(t) -a Z^2(t))dt + G Zd\mathbb{W}(t). \label{eq1:lv}
\end{equation}        
where $a$, $b$ and $G$ are constants. The general solution of Eq.(\ref{eq1:lv}) is given by \cite{arat2003,gard}
\begin{eqnarray}
Z(t) &=& \exp \left \{ \left ( b - {G^2\over 2}\right ) t +G\mathbb{W}(t) \right \} \nonumber \\
& & \times \left [ {1\over Z(0)} + a  \int_0^t \exp \left \{ \left ( b - {G^2\over 2}   \right )s + G \mathbb{W}(s) \right \} ds \right ]^{-1} .
\label{eq2:lv}
\end{eqnarray}
By treating $Z$ as a hypecomplex-valued stochastic process, we may cast Eq.(\ref{eq1:lv}) into a system of stochastic Lotka-Volterra whose solution is obtain by projecting the solution (\ref{eq2:lv}) along the units on the underlying hypercomplex. Indeed, set $Z = \sum_{i=1}^n X_i \,e_i$, $a = \sum_{i=1}^n a_i\, e_i$, $b =\sum_{i=1}^n b_i\, e_i$, $G = \sum_{i=1}^n G_i\, e_i$,  and $\mathbb{W} = \sum_{i=1}^n W_i\, e_i$, where the $e_i$'s are units of an $n$-dimensional commutative hypercomplex $H$ and the $W_i$'s are independent standard Wiener processes. Then, Eq. (\ref{eq1:lv}) is equivalent to 
\begin{eqnarray}
dX_i(t) &=& \left (\sum_{j,k =1}^n \gamma_{kji}b_k X_j -\sum_{j,k,r,s=1}^n \gamma_{rjk}\gamma_{ski}b_sX_rX_j \right ) dt
\nonumber \\ 
& & + \sum_{j,k,r,s=1}^n \gamma_{sjk}\gamma_{kri} G_s X_jdW_r(t),\; i= 1: n. \label{eq3:lv}     
\end{eqnarray} 
Thus, we have proved the following result.
\begin{theo}
The  stochastic Lotka-Volterra system (\ref{eq3:lv}) is integrable by qua- \\ dradures provided the  
$\gamma_{ijk}$'s appearing in it are structure constants of a commutative hypercomplex. Its solution is then obtained by projecting Eq. (\ref{eq2:lv}) along the units of the underlying commutative hypercomplex.  
\end{theo}
In guise of illustration, we take $H = \mathbb{C}_p$. Then, the system (\ref{eq3:lv}) assumes the form
\begin{eqnarray}
 dX_1(t) &=& [b_1X_1(t) + p\, b_2 X_2(t) - 2\,p\,a_2 X_1(t)X_2(t) -a_1 (X_1^2(t) + p X_2^2(t))]dt \nonumber \\
 & & + [G_1 X_1(t) + p\, G_2 X_2(t)]dW_1(t) \nonumber \\
 & & + p\,[G_2 X_1(t) + G_1 X_1(t)]dW_2(t), \label{eq4:lv}\\
 dX_2(t) & =&  [b_2 X_1(t) + b_1 X_2(t) - 2\,a_1\,X_1(t) X_2(t) - a_2 (X_1^2(t) + p X_2^2(t))]dt \nonumber \\
& & + [G_2 X_1(t) + G_1 X_2(t)]dW_1(t) \nonumber \\
& & + [G_1 X_1(t) + p\, G_2 X_2(t)]dW_2(t). \label{eq5:lv} 
\end{eqnarray} 
Through  extracting the real and imaginary, we infer from Eq. (\ref{eq2:lv}) after some calculations, that the solution of the system (\ref{eq4:lv})-(\ref{eq5:lv}) takes the form 
\begin{eqnarray}
X_1(t) & = & {\mbox{e}^{\alpha(t)} \over \gamma^2(t) -p\, \delta^2(t)}\;
\left \{\gamma(t) \cos_p[\beta(t)] - p\, \delta(t) \sin_p[\beta(t)]\right\}
, \label{eq6:lv}\\
X_2(t) & =& {\mbox{e}^{\alpha(t)} \over \gamma^2(t) -p\, \delta^2(t)}\;
\left \{\gamma(t) \sin_p[\beta(t)] - \, \delta(t) \cos_p[\beta(t)]\right \} \label{eq7:lv}
\end{eqnarray}  
where the functions $\alpha$, $\beta$, $\gamma$, and $\delta$ are given by
\begin{eqnarray}
\alpha(t) &=& \left ( b_1 -{G_1^2\over 2} - {p\,G_2^2\over 2}   \right ) t +
G_1\, W_1(t) + p\, G_2\, W_2(t), \label{eq8:lv}\\
\beta(t) &=& (b_2 -G_1\,G_2)t + G_2 W_1(t) + G_1\, W_2(t), \label{eq9:lv}\\
\gamma(t) & = & {X_1(0)\over X_1^2(0)-pX_2^2(0)} + a_1 \int_0^t \mbox{e}^{\alpha(s)}
\cos_p[\beta(s)]\,ds \nonumber \\
& &  + p\, a_2 \int_0^t \mbox{e}^{\alpha(s)}\sin_p [\beta(s)]\,ds , \label{eq10:lv}\\
\delta(t) & = & -{X_2(0)\over X_1^2(0) - p X_2^2(0)}
 + a_1 \int_0^t \mbox{e}^{\alpha(s)}
\sin_p[\beta(s)]\,ds \nonumber \\
& &  +  a_2 \int_0^t \mbox{e}^{\alpha(s)}\cos_p [\beta(s)]\,ds . \label{eq11:lv}
\end{eqnarray}  

\section{Conclusion}
In this work we have demonstrated that we can exactly solve some systems of SDEs through hypercomplexification provided the base equation is solvable. We have explicitly treated few examples to showcase the efficiency of the procedure. Precisely, we exhibited integrable systems of linear SDEs, obtained linearization criteria for some systems of SDEs and we showed how integrable stochastic Lotka-Volterra can be begotten. We stress that although we have limited ourselves to few cases, hypercomplexification applies to other types of SDEs not considered here. In particular, the method can be applied to SDEs driven by non-white noises. The diligent reader may have noticed that hypercomplexification is fundamentally an algebraic method which is enabled by the analyticity of data. Computer algebra packages for finding  exact solutions of  scalar SDEs may be extend to handle systems of SDEs that can be hypercomplexified to a solvable base equation. Such an endeavor may be undertaken at the cost of very little programming effort. 

There are few problems we could not solve that deserve a future attention. Some systems of SDEs may be hypercomplexifiable only after an appropriate change of variables. It will be interesting to characterize such systems of SDEs. Also, we have established that the linear system (\ref{eq10:lsde}) is integrable provided the $\gamma_{ij}^k$ are structure constants of a commutative hypercomplex. It would be interesting to find the necessary and sufficient conditions under which a system of SDEs $$dY_i(t) = a_i(t, Y_{1:n}(t)) dt + \sum_{j=1}^n b_{ij}(t, Y_{1:n}(t))dW_j(t),\; i = 1:n,$$ can be invertibly mapped to the integrable linear system  (\ref{eq10:lsde}) through a change of variables $Y_i = h_i(t, Y_{1:n}(t)),\; i=1:n$.
Even the case $n= 2$ leads to  formidable computations that we hope can be handled by an appropriate computer algebra package. 

\ack
This work  was completed during the visit of CWS at the DST-NRF Centre of Excellence in Mathematical and Statistical Sciences of the University of the Witwatersarand, Johannesburg, South Africa. The authors thankfully acknowledge financial support from the NRF of South Africa.
\section*{References}
\bibliography{sde}
\end{document}